\documentclass[12pt]{article}
\usepackage{graphicx}  
\usepackage{amsmath} 
\usepackage{amssymb} 
\usepackage{amsfonts} 
\usepackage{amsthm} 
\usepackage{longtable} 
\usepackage{multirow}
\usepackage{color}
\input{epsf}
\input xy
\xyoption{all}
\usepackage{latexsym}
\newtheorem{lemma}{Lemma}[section] 
\newtheorem{theorem}[lemma]{Theorem} 
\newtheorem{proposition}[lemma]{Proposition}

\newtheorem{definition}[lemma]{Definition} 
\theoremstyle{remark}
\newtheorem{remark}[lemma]{Remark}

\renewcommand{\Gamma}{\varGamma} 
\renewcommand{\epsilon}{\varepsilon} 
\renewcommand{\bar}{\overline} 
\renewcommand{\hat}{\widehat} 

\renewcommand{\leq}{\leqslant}

\newcommand{\Po}{\mathcal{P}}
\newcommand{\G}{\mathcal{G}}
\newcommand{\M}{\mathcal{M}}
\newcommand{\R}{\mathcal{R}}
\newcommand{\Lm}{\mathcal{L}}
\newcommand{\C}{\mathcal{C}}
\newcommand{\Oh}{\mathcal{O}}
\newcommand{\Q}{\mathcal{Q}} 
 
\newcommand{\T}{\mathcal{T}}

\newcommand{\ep}{\varepsilon}

\begin{document}

\title{On Roli's Cube}
 
\author{Barry Monson\thanks{This work was generously supported at the Universidad Nacional Aut\'{o}noma de M\'{e}xico (Morelia) by
PAPIIT-UNAM  grant \#IN100518 .} \\ 
University of New Brunswick\\
Fredericton, New Brunswick, Canada E3B 5A3}

\date{ \today }
\maketitle  

\begin{abstract}
First described in 2014, \textit{Roli's cube}  $\R$ is a 
chiral $4$-polytope, faithfully realized in Euclidean $4$-space (a situation  earlier
thought to  be impossible).  Here we describe $\R$ in a new way,  determine  its minimal
regular cover,  and reveal connections to
the M\"{o}bius-Kantor configuration.



\end{abstract}

\medskip

\section{Introduction}\label{intro}

Actually \textit{Roli's cube} $\R$ isn't a cube, although it does share the $1$-skeleton of a $4$-cube. First
described by Javier (Roli) Bracho, Isabel Hubard and Daniel Pellicer in \cite{bracho:2014aa}, 
$\R$ is a chiral $4$-polytope of type $\{8,3,3\}$, faithfully realized in $\mathbb{E}^4$
(a situation earlier thought impossible).
Of course, Roli didn't himself name $\R$; but  the  eponym is pleasing to his colleagues and has taken hold.

Chiral polytopes  with realizations of `full rank' had  (incorrectly) 
been shown  not to exist by
Peter McMullen in \cite[Theorem 11.2]{fullrank}. Mind you,  these objects 
do seem to be elusive. 
Pellicer has proved in \cite{chir4and5}  that chiral polytopes of full rank  can exist only in ranks
$4$ or $5$.

Roli's cube $\R$ was  constructed in  \cite{bracho:2014aa} as a \textit{colourful polytope}, starting from 
a hemi-$4$-cube in projective $3$-space. (For more on this, see  Section~\ref{copo}.)
The construction given here in Section~\ref{rolq}
is a bit different, though
certainly  closely related. In Section~\ref{mincov} we can then easily  manufacture the minimal
regular cover $\T$ for $\R$, and give both a presentation and faithful representation for
its automorphism group. Along the way, we encounter both the 
M\"{o}bius-Kantor Configuration $8_3$
and the regular complex polygon   $3\{3\}3$.

In what follows, we can make our way with  concrete examples, so we won't need much
of the general theory of abstract regular or chiral polytopes
and their realizations. We refer the reader to \cite{arp}, \cite{GeRP} and \cite{schulte1} for  more.

\section{The $4$-cube: convex, abstract and colourful}\label{4cub} 

\begin{figure}[ht]\begin{center}
\includegraphics[width=50mm]{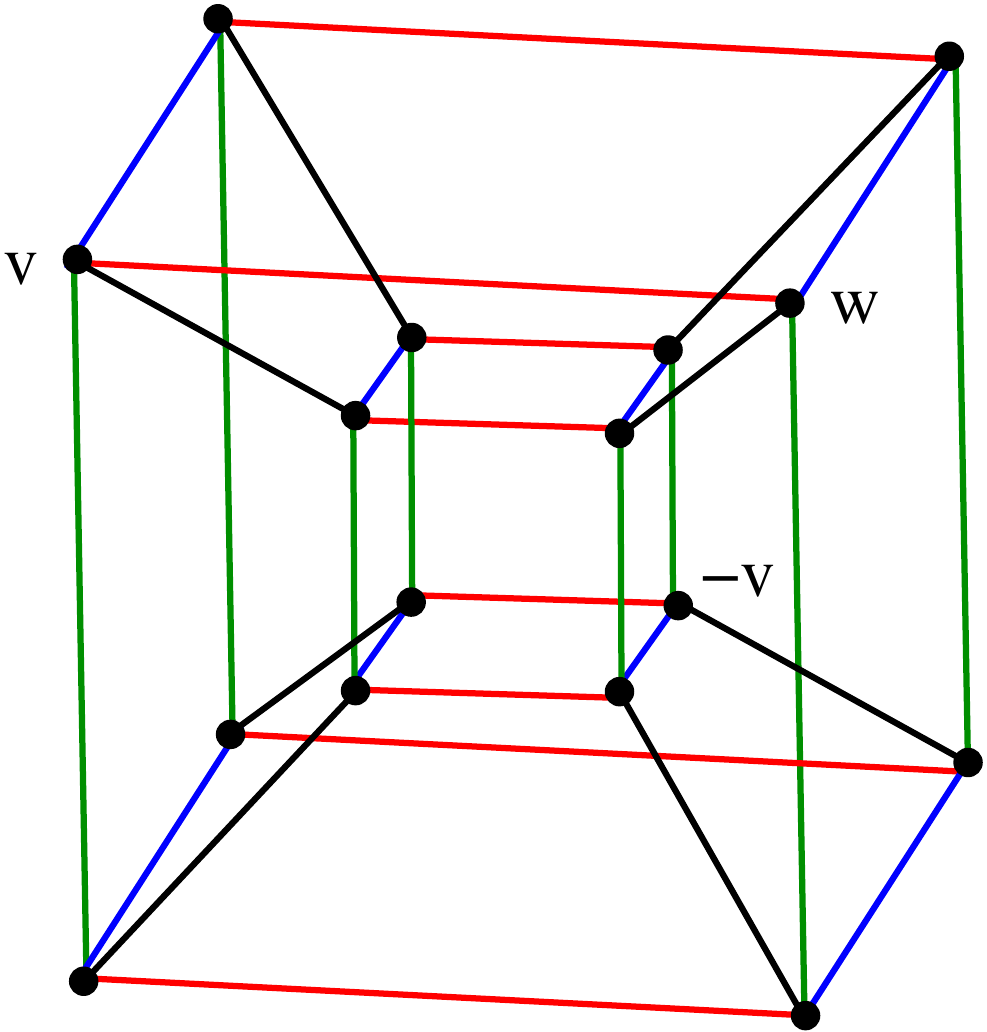}
\end{center}
\caption{A $2$-dimensional look at a  $3$-dimensional projection of the $4$-cube.}
\label{proj2}
\end{figure}

The most familiar of the regular convex polytopes in Euclidean space 
$\mathbb{E}^4$ is surely the $4$-cube $\Po= \{4,3,3\}$. 
A  familiar projection of $\Po$ into $\mathbb{E}^3$ is displayed in 
Figure~\ref{proj2}.

\begin{figure}[ht]\begin{center}
\includegraphics[width=120mm]{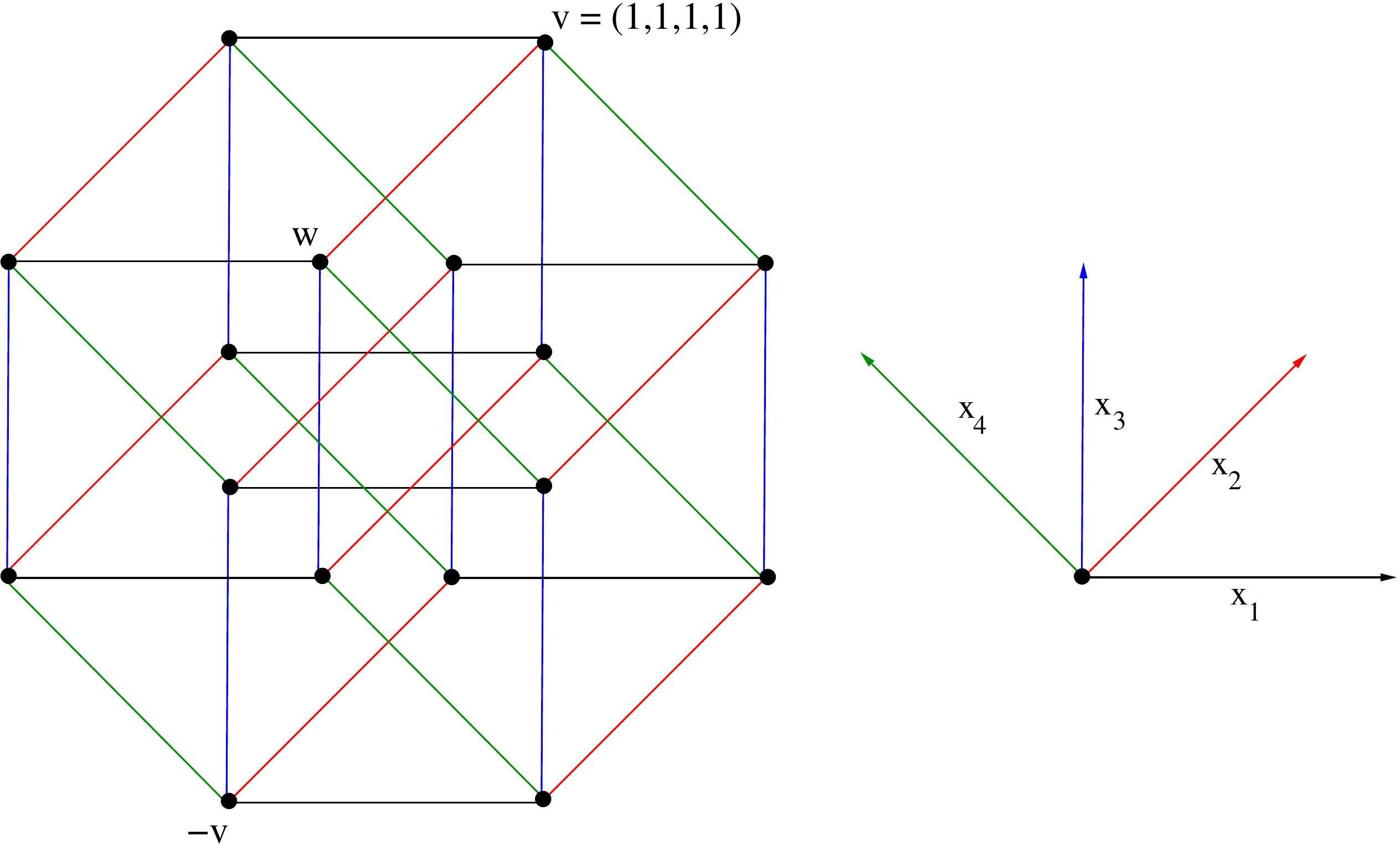}
\end{center}
\caption{The most symmetric $2$-dimensional projection of the $4$-cube.}
\label{proj1}
\end{figure}

Let us equip $\mathbb{E}^4$ with its usual basis $b_1,\ldots,b_4$ and inner product.
Then we may take the vertices of $\Po$ to be the $16$ sign change vectors
 \begin{equation}\label{vert} e = (\ep_1, \ldots,\ep_4) \in \{\pm 1\}^4.\end{equation}
 At any such vertex there is an edge (of length $2$) running 
in each of the $4$ 
coordinate directions, so that  $\Po$ has $32 = \frac{16\cdot 4}{2}$  edges. 
Similarly we count the 
$24$  squares $\{4\}$  as faces of  dimension $2$.
Finally, $\Po$ has $8$ facets; these faces of   dimension $3$  
are ordinary cubes  $\{4,3\}$. They lie in four pairs
of supporting hyperplanes   
orthogonal to the coordinate axes. It is enjoyable to hunt for these faces in Figure~\ref{proj1}, where the $8$ 
parallel edges in each of the coordinate directions have colours black, red, blue and green, respectively.
 
%
%

Let us turn to the symmetry group $G$ for $\Po$. Each symmetry $\gamma$ is determined by its action on the
vertices, which clearly can be permuted with sign changes in all possible ways. Thus $G$ has order
$2^4\cdot 4! = 384$, and we may think of it  as being comprised of all $4 \times 4$ signed permutation matrices.

In fact, $G$ can be generated by reflections $\rho_0, \rho_1, \rho_2, \rho_{3}$ in hyperplanes. Here
$\rho_0$ negates the first coordinate $x_1$ (reflection in the coordinate hyperplane  orthogonal to $b_1$);
and,  for $1 \leq j \leq 3$, $\rho_j$ transposes coordinates $x_j, x_{j+1}$ (reflection  in the hyperlane  orthogonal
to  $b_j-b_{j+1}$).  

Note that the reflection in the $j$-th coordinate  hyperplane is 
$  \rho_0^{\rho_1 \cdots \rho_{j-1}} $ for $1\leq j\leq 4$. 
(We use the notation $\gamma^{\eta} := \eta^{-1} \gamma \eta$.)
The product of these $4$ special reflections, in any order,
 is the central element $\zeta: t \mapsto -t $. It is easy to check as well that
\begin{equation}\label{cent}
\zeta = (\rho_0\rho_1 \rho_2\rho_{3})^4\;.
\end{equation}
The \textit{Petrie symmetry} $\pi = \rho_0\rho_1 \rho_2\rho_{3}$ therefore has period $8$.

For purposes of  calculation, we note that 
$G \simeq C_2^4 \rtimes S_4$ is   a semidirect product. Under this isomorphism,
 each  $\gamma \in G$ factors uniquely as 
$\gamma =  e \mu $, where $\mu  \in S_n$ is a permutation of  
$\{1,\ldots,4\}$ (labelling the  coordinates);
and $e$ is a sign change vector, as in (\ref{vert}). 
Note that 
$$ e^{\mu} =  (\ep_1, \ep_2, \ep_3, \ep_4)^{\mu} = (\ep_{(1)\mu^{-1}}, \ep_{(2)\mu^{-1}}, 
\ep_{(3)\mu^{-1}} , \ep_{(4)\mu^{-1}})\, .$$
%
Now really $\gamma$ is a signed permutation matrix. But it is convenient to abuse notation, 
keeping in mind that each $e$ corresponds to a diagonal matrix of signs
and each $\mu$ to a permutation matrix. Thus we might write
\begin{equation}\label{pet1}
\pi   = \rho_0\rho_1 \rho_2\rho_{3} = ( -1,1, 1,1)\cdot (4,3,2,1) =\left[
\begin{array}{rrrr}
 0&0&0&-1 \\
 1&0&0&0 \\
 0&1&0&0 \\
 0&0&1&0
\end{array}
\right]\;.
\end{equation}

%
Next we use the group $G =  \langle \rho_0, \rho_1, \rho_2, \rho_{3}\rangle$ to remanufacture the cube. 
In this (geometric) version
of \textit{Wythoff's construction} \cite[\S 2.4]{rcp} we  choose a \textit{base vertex} $v$ fixed by the
subgroup  $G_0 := \langle \rho_1, \rho_2, \rho_{3}\rangle$ (which permutes the coordinates in all ways). Thus,
$v =  c(1,1,1,1)$ for some $c \in \mathbb{R}$. To avoid a trivial construction we take $c \neq 0$, so, up to 
similarity, we may
use $c=1$. Then the orbit of $v$ under $G$ is just the set of $16$ points in (\ref{vert}); and their convex hull
returns $\Po$ to us. Since $G_0$ is the full stabilizer of $v$ in $G$, the vertices correspond to right 
cosets $G_0\gamma$.

The beauty of Wythoff's construction is that all faces of  $\Po$   can be constructed in a similar way by induction on 
dimension (\cite[Section 1B]{arp}, \cite{coxeter10} and \cite{GeRP}). For example, the vertices $v = (1,1,1,1)$ and $v\rho_0 = (-1,1,1,1)$
of the \textit{base edge} of $\Po$ 
are just the orbit of $v$ under the subgroup $G_1 := \langle \rho_0, \rho_2, \rho_{3}\rangle$; and edges
of $\Po$ correspond to right cosets of the new subgroup $G_1$. Furthermore, a more careful look reveals that
a vertex is incident with an edge just when the  corresponding cosets have non-trivial intersection.

Pursuing this, we see that the face lattice of $\Po$ can be recontructed as a coset geometry based 
on   subgroups 

\begin{equation}\label{dist}
G_0, G_1, G_2, G_3,  \;\;\mathrm{where}\;G_j := \langle \rho_0,\ldots, \hat{\rho_{j}},\ldots, \rho_{3}\rangle.
\end{equation}
From this point of view, $\Po$ becomes an 
\textit{abstract regular $4$-polytope}, a partially ordered set whose automorphism group is  $G$.
Notice that the distinguished subgroups in (\ref{dist}) provide the 
proper faces in a \textit{flag} in $\Po$, namely a mutually incident vertex, edge, square and $3$-cube.

The crucial structural property of $G$ is that it should be 
a string C-group with respect to the generators $\rho_j$. 
A \textit{string C-group}  is a quotient 
of a  Coxeter group with linear diagram  under which an `intersection condition' on   subgroups
generated by subsets of generators,  such as  those in (\ref{dist}),  
is preserved  \cite[Sections 2E]{arp}. 

For the $4$-cube $\Po$, $G$ is actually isomorphic to the 
Coxeter group $B_4$ with diagram 

\begin{equation}\label{B4diag}{\bullet}\frac{4}{\;\;\;\;\;\;\;} {\bullet}\frac{3}{\;\;\;\;\;\;\;}{\bullet}\frac{3}{\;\;\;\;\;\;\;}
{\bullet}
\end{equation} 
Comparing the geometric and abstract points of view, we say that the convex $4$-cube 
is a \textit{realization} of its face lattice (the abstract $4$-cube).

When we think of a polytope from the abstract point of view, we often use the term \textit{rank} instead
of `dimension'.   An abstract polytope $\Q$ is said to be \textit{regular} if its automorphism group 
is transitive on flags (maximal chains in $\Q$).
 Intuitively, regular polytopes have maximal symmetry (by reflections). Next up
 are \textit{chiral} polytopes, with exactly two flag orbits and such that adjacent flags are
 always in different orbits (so maximal symmetry by rotations, but without reflections).

We will  soon encounter  less familiar abstract regular or chiral polytopes, with their
realizations. %
For a first example, suppose that we map  (by central projection)
the  faces of $\Po$ onto the  $3$-sphere $\mathbb{S}^3$
centred at the origin. We can then reinterpret $\Po$ as a regular \textit{spherical polytope} (or tessellation),
with the same symmetry group $G$. Now the centre of $G$ is 
the subgroup $\langle \zeta\rangle$ of order $2$. The quotient group $G/\langle \zeta\rangle$
has order $192$ and is still a string C-group. The corresponding regular polytope
is the \textit{hemi-$4$-cube} $\mathcal{H} = \{4,3,3\}_4$,  now realized in projective
space $\mathbb{P}^3$ \cite[Section 6C]{arp}; see Figure~\ref{k44fig}. 
By (\ref{cent}), the product of the four
generators of $G/\langle \zeta\rangle$ has order $4$; this is recorded as the  subscript in the
Schl\"{a}fli symbol for $\mathcal{H}$.

Now we can   outline the construction of
Roli's cube given in \cite{bracho:2014aa}.

\medskip
\section{Colourful polyopes}\label{copo}
The image in Figure~\ref{proj2} or on the left in Figure~\ref{proj1} can just as well be understood as a graph $\G$,
namely the \textit{$1$-skeleton} of the $4$-cube $\Po$. In fact, we can recreate the 
abstract (or  combinatorial) structure of $\Po$ from just the edge colouring of $\G$: 
for $0\leq j \leq 4$, the $j$-faces of $\Po$ can be identified with the components of those subgraphs
obtained by keeping  just  edges with some selection of the $j$ colours (over all such choices). 
We  therefore say that  $\Po$ is a \textit{colourful polytope}.

Such polytopes were introduced in
\cite{schulteGID3}.  In general, one begins with a finite, connected $d$-valent graph $\G$
admitting a (proper) edge colouring, say by the symbols $1,\ldots,d$.
Thus each of the colours provides a $1$-factor for $\G$.
The graph  $\G$ determines an (abstract)  colourful polytope $\Po_{\G}$ as follows. For $0\leq j \leq d$, 
a typical $j$-face $(C,v)$ is identified with the set  all vertices of $\G$ 
connected to a given vertex $v$ by a path using
only colours   from some subset $C$ of size $j$ taken from $\{1,\ldots,d\}$. The $j$-face $(C,v)$
is incident with the $k$-face $(D,w)$ just when $C \subseteq D$ and $w$ can be reached from $v$
by a $D$-coloured path. (This means  that $j\leq k$;  and we can just as well take $w=v$. 
The minimal face of  rank $-1$ in $\Po_{\G}$  is formal.)  Notice that $\Po_{\G}$ is a \textit{simple}
$d$-polytope whose $1$-skeleton   is just $\G$ itself.
From \cite[Theorem 4.1]{schulteGID3}, the automorphism group of $\Po_{\G}$ is isomorphic to the group
of colour-preserving graph automorphisms of $\G$. (Such automorphisms are allowed to 
permute the $1$-factors.)

\begin{figure}[ht]\begin{center}
\includegraphics[height=60mm]{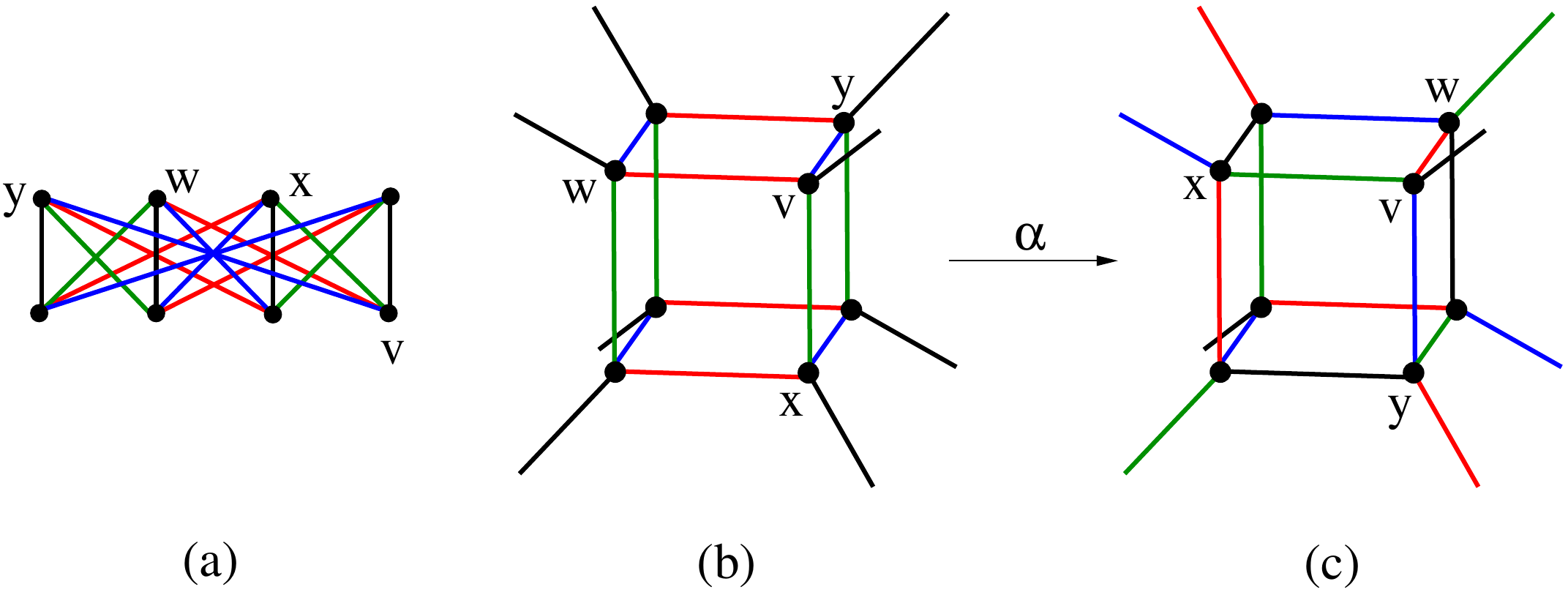}
\end{center}
\caption{The   graph $K_{4,4}$ in (a)  is the $1$-skeleton of the hemi-$4$-cube $\{4,3,3\}_4$ (b),(c).}
\label{k44fig}
\end{figure}

It is easy to see that the hemi-$4$-cube $\mathcal{H}$ is also colourful. Its $1$-skeleton is the 
complete bipartite graph $K_{4,4}$ found in  
Figure~\ref{k44fig}.  
We   obtain this  graph from  Figure~\ref{proj2} or Figure~\ref{proj1} by   identifying antipodal pairs of points,
like $v$ and $-v$. 

If we lift $\mathcal{H}$, as it is now, to $\mathbb{S}^3$, we regain the coloured
$4$-cube $\Po$.  Now  keep $K_{4,4}$ embedded in $\mathbb{P}^3$, as in Figure~\ref{k44fig}.
But,  following    \cite{bracho:2014aa}, observe that $K_{4,4}$ admits the automorphism $\alpha$
which cyclically permutes, say, the first three vertices $y,w,x$ in the top block, leaving the rest fixed.
Clearly, $\alpha$ is a non-colour-preserving automorphism of $K_{4,4}$, so its effect
is to recolour $12$ of the edges in the embedded graph.  
On the abstract level nothing has 
changed for the resulting colourful  polytope; it is still the hemi-$4$-cube  $\mathcal{H}$. 
But faces of ranks $2$ and $3$
are now differently embedded in $\mathbb{P}^3$.
For example, the red-blue $2$-face on $v$, which  is planar in Figure~\ref{k44fig}(b),   becomes
 a helical quadrangle Figure~\ref{k44fig}(c) and   thereby acquires an orientation.
 According to Definition~\ref{petdef}, these helical polygons are Petrie
 polygons for the standard realization of $\mathcal{H}$ in Figure~\ref{k44fig}(b).

The newly coloured geometric object, which we might label $\mathcal{H}^R$, is   a 
\textit{chiral realization}  of the abstract regular polytope $\mathcal{H}$.   Comforted by the fact that $\mathbb{P}^3$ is orientable,  we could just as well apply
$\alpha^{-1}$ to obtain the left-handed version $\mathcal{H}^L$. These two \textit{enantiomorphs}
are oppositely embedded in $\mathbb{P}^3$, though both remain isomorphic to $\mathcal{H}$ 
as partially ordered sets.
If we lift either enantiomorph to $\mathbb{S}^3$, we obtain a chiral
$4$-polytope faithfully realized in $\mathbb{E}^4$ \cite[Theorem 2]{bracho:2014aa}. This is Roli's cube $\R$.

Next we set the stage for a slightly different construction of $\R$, 
without the use of $\mathbb{P}^3$.

\section{Petrie polygons of the $4$-cube}\label{petri}
Let us consider the progress of the base vertex $v = (1,1,1,1)$ as we apply successive powers of $\pi$ in (\ref{pet1}).
We get a centrally symmetric $8$-cycle of vertices
$$v\rightarrow(1,1,1,-1)\rightarrow (1,1,-1,-1)\rightarrow (1,-1,-1,-1)\rightarrow
 -v = (-1,-1,-1,-1)\rightarrow \ldots  \;.
$$
Starting from $v$ in Figure~\ref{proj1} we therefore proceed in   coordinate directions 4, 3, 2, 1
(indicated by different colours), then repeat again.  This traces   out the peripheral octagon $\C$,
which in fact is a Petrie polygon for $\Po$.

\begin{definition}\label{petdef}
A \textit{Petrie polygon} of a $3$-polytope is an edge-path such that
any $2$ consecutive edges, but no $3$, belong to a $2$-face. 
We then say that   a \textit{Petrie polygon} of a $4$-polytope  $\Q$ is   an edge-path such that
any $3$ consecutive edges, but no $4$, belong to \emph{(}a Petrie polygon of\emph{)} a facet of $\Q$. 
\end{definition}
For the cube $\Po$, the parenthetical condition is actually superfluous; compare  \cite{coxweissA}.

Clearly, we can begin a Petrie polygon at any vertex, taking any of the $4!$ orderings of the colours. But this counts
each octagon in $16$ ways. We conclude that $\Po$ has $24$ Petrie polygons. What we really use here is
the fact that $G$ is transitive on vertices, and that at any fixed vertex, $G$ permutes the edges in all 
possible ways. We see that $G$ acts transitively on Petrie polygons.

But the (global) stabilizer of $\mathcal{C}$ (constructed above with the help of  $\pi$ and $v$) is the dihedral group $K$ of  order $16$  generated by 
$\mu_0  = \rho_0\rho_2\rho_3\rho_2 = (-1,1,1,1)\cdot(2,4),$ and 
$\mu_1 = \mu_0\pi = \rho_2\rho_3\rho_2\rho_1\rho_2\rho_3  = (1,1,1,1)\cdot (1,4)(2,3)$.
(Such  calculations are routine using either signed permutation matrices or the decomposition in
$C_2^4 \rtimes S_4$. Note that any $4$ consecutive vertices of $\mathcal{C}$ form a 
basis of $\mathbb{E}^4$.) We confirm that $\Po$ has $24 =384/16$ Petrie polygons.

Now we move to  the \textit{rotation subgroup} 
$$G^+ = \langle \rho_0\rho_1, \rho_1\rho_2,  \rho_2\rho_3 \rangle.$$
It has order $192$ and consists of the signed permutation matrices of determinant $+1$.
Note that $K< G^+$. Thus, under the action of $G^+$, there are two orbits of Petrie 
polygons of $12$ each. Let's label these two \textit{chiral classes}
$R$ and $L$ for right- and left-handed, taking $\mathcal{C}$ in class $R$.

The two chiral classes must be swapped
by any non-rotation, such as any  $\rho_j$. To distinguish
them,  we could take the determinant of the matrix whose rows are any $4$ consecutive vertices on 
a Petrie polygon. The two chiral classes $R$ and $L$ then have determinants $+8$, respectively, $-8$.
Or starting from a common vertex, the edge-colour sequence along a polygon in one class is an 
odd permutation of the colour sequence for a polygon in the other class.

The inner octagram  ${\C}^*$ in Figure~\ref{proj1} is another Petrie polygon. Start at the vertex 
$w = (v)\rho_1\rho_0\rho_1 = (1,-1,1,1)$ 
which is adjacent to $v$ along a red edge; then proceed in directions 4, 1, 2, 3 and repeat. However, the remaining
Petrie polygons appear in less symmetrical fashion in  Figure~\ref{proj1}.

Note that $\mu_0$ actually acts on the diagram in Figure~\ref{proj1} as a reflection in a vertical line,
whereas $\pi$ rotates the octagon $\mathcal{C}$ and octagram $\mathcal{C}^*$  in opposite senses. On the other hand, $\mu_2 =  \rho_1\rho_2\rho_0\rho_1 = (1,-1,1,1)\cdot (1,3)$ is an element of $G^+$
which swaps $\mathcal{C}$ and  $\mathcal{C}^*$. Thus there are $6$ such unordered pairs
like $\mathcal{C}, \mathcal{C}^*$ in class $R$ and another $6$  pairs in  class $L$.

\medskip
\begin{remark}
It can be shown that  Figure~\ref{proj1} is the most symmetric orthogonal projection of $\Po$ to a plane
\cite[\S 13.3]{rp}. Since all edges after projection have a common length, we may say that this projection
is \textit{isometric}. 
 
The Petrie symmetry $\pi$ is one instance of a  \textit{Coxeter element}  in the group $G = B_4$, namely 
a product   of the four generators in some order. 
 All such Coxeter elements are conjugate.  Each of them has invariant planes
 which give rise to the sort of orthogonal projection displayed in Figure~\ref{proj1}.
  A procedure for finding these planes is detailed in \cite[3.17]{humph}.
For $\pi$,     the
two   planes are spanned by the rows of 
$$\left[\begin{array}{rrrr}
 \frac{1}{\sqrt{2}}&  \frac{1}{2} & 0&  \frac{-1}{2} \\
0&  \frac{1}{2} &  \frac{1}{\sqrt{2}} & \frac{1}{2}   \end{array}\right]
\;\mathrm{and}\;
\left[\begin{array}{rrrr}
 \frac{1}{\sqrt{2}}&  \frac{-1}{2} & 0&  \frac{1}{2} \\
0&  \frac{1}{2} &  \frac{-1}{\sqrt{2}} & \frac{1}{2}   \end{array}\right] 
.$$
These   planes are  orthogonal complements; and $\pi$ acts on them by 
rotations through $45^{\circ}$ and $135^{\circ}$, respectively. Figure~\ref{proj1} results from 
projecting $\Po$ onto the first plane.
\hfill$\square$
\end{remark}

\bigskip
\section{The map $\mathcal{M}$ and the M\"{o}bius-Kantor Configuration $8_3$}\label{83con}

\begin{figure}[ht]\begin{center}
\includegraphics[width=70mm]{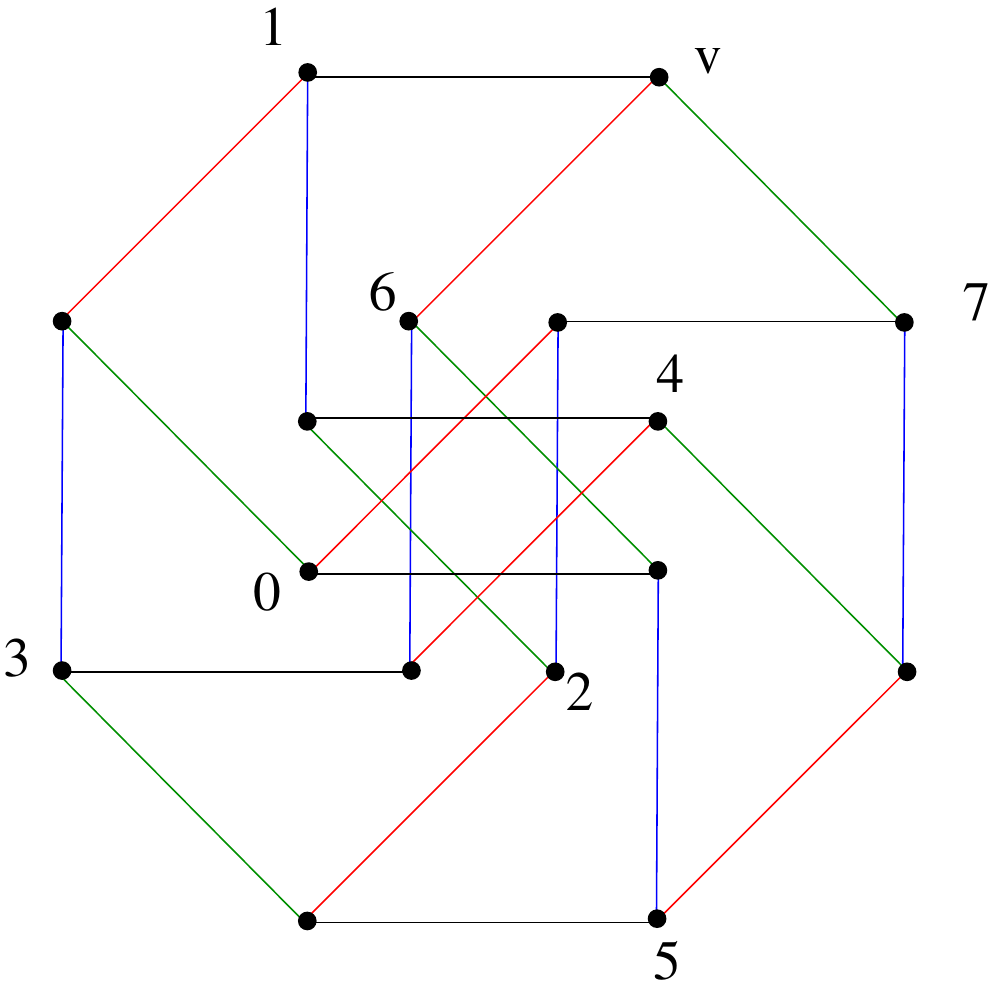}
\end{center}
\caption{The Levi graph $\Lm = \{8\} + \{\frac{8}{3}\}$ for the configuration $8_3$.}
\label{levi}
\end{figure}

Look again at the companion Petrie polygons
$\C,\C^*$ in Figure~\ref{proj1}.
Now working around the rim clockwise from $v$ delete the edges coloured
blue, red, black, green, and repeat. We are left with the trivalent graph $\Lm$
displayed in Figure~\ref{levi}.

In fact, $\Lm$ is the \textit{generalized Petersen graph} $\{8\} + \{\frac{8}{3}\}$, studied in detail 
by Coxeter in \cite[Section 5]{Coxeter:1950aa}. The graph is $2$-arc transitive, so that its automorphism group
has order $96 = 16\cdot 3\cdot 2^{2-1}$ \cite[Chapter 18]{biggsAGT}. We return to this group later.

We have labelled alternate vertices of $\Lm$ by the residues $0,1,2,3,4,5,6,7\pmod{8}$.
These will  represent the points in a M\"{o}bius-Kantor configuration $8_3$. The remaining 
(unlabelled) vertices of $\Lm$ represent the $8$ lines in the configuration. Thus we have lines
$013$ (represented by the `north-west' vertex $(-1,-1,1,1)$), $124$, $235$, and so on, including
line $671$ represented by $v$.

Notice that we can interpret the configuration as being comprised of two quadrangles with vertices $0,2,4,6$ and $1,3,5,7$,
\textit{each inscribed in the other}: vertex $0$ lies on edge $13$, vertex $1$ lies on edge $24$, and so on.

So far this configuration $8_3$ is purely abstract. In fact, it can be realized as
a point-line configuration in a projective (or affine) plane 
over any field in which $$z^2 -z +1 = 0$$
has a root, certainly over $\mathbb{C}$. However, $8_3$  cannot be realized
in the real  plane.

Coxeter made other observations in \cite{Coxeter:1950aa}, including the fact  that the graph $\Lm$ is
a sub-$1$-skeleton of the $4$-cube. Altogether $\Lm$ contains 
 $6$ Petrie polygons, which we can briefly describe by their alternate vertices:
$$
\begin{array}{ccc}
 0246 (= \mathcal{C}^*)\;\; &0541 & 1256 \\
 1357 (=\mathcal{C})\;\; &2367 & 0743 \\
\end{array}
$$
Hence,  the configuration can be regarded as a pair of mutually inscribed quadrangles
in  three  ways.

Observe   that each edge of $\Lm$ lies on exactly two of the $6$ octagons. For example, the top edge
with vertices labelled $1$ and $v$ lies on octagons $1357$ and $1256$. (It does not matter that two such
octagons then share a second edge opposite the first.)
Furthermore, each vertex lies on the three octagons determined by
choices of two edges.  We can thereby construct a $3$-polytope $\M$ of type $\{8,3\}$, with $6$ octagonal 
faces, whose $1$-skeleton is $\Lm$.  In short, $\M$ is realized by substructures of the $4$-cube $\Po$.

\begin{figure}[ht]\begin{center}
\includegraphics[width=90mm]{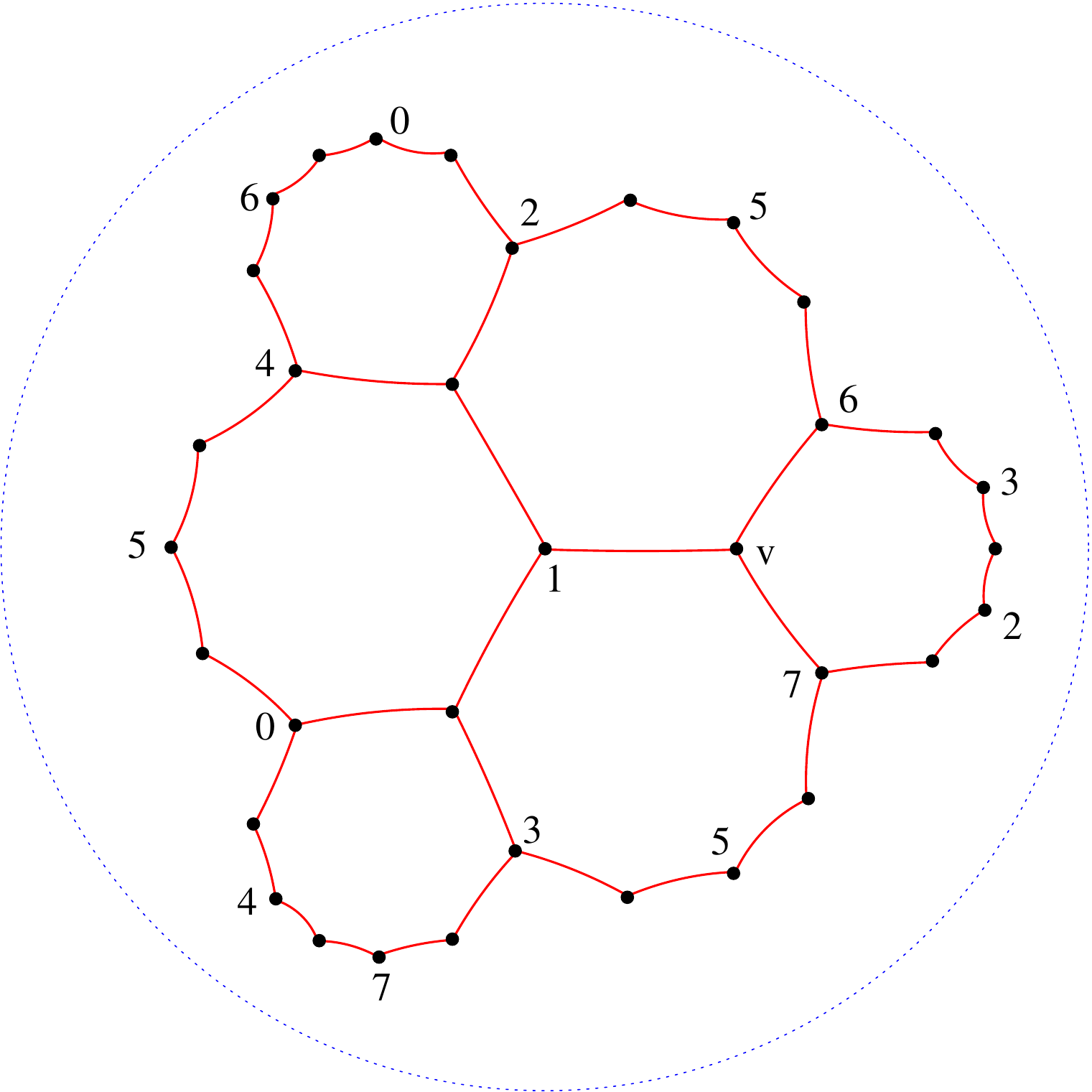}
\end{center}
\caption{Part of the tessellation $\{8,3\}$ of the hyperbolic plane.}
\label{83map}
\end{figure}

Moving sideways, we can reinterpret $\M$ in a more familiar topological  way as a map
on a compact  orientable surface of genus $2$. Recall that  $\M$ is  covered by the tessellation $\{8,3\}$ of the hyperbolic plane, as indicated 
in Figure~\ref{83map}.

Now return to $\mathbb{E}^4$ where the combinatorial structure of $\M$ is handed  to us as faithfully realized. 
Drawing on  \cite[Section 8.1]{coxmos}, we have that the rotation group $\Gamma(\M)^+$ for $\M$
is generated by two special Euclidean symmetries:
 
$\sigma_1 = \pi = \rho_0\rho_1\rho_2\rho_3  = (-1,1,1,1)\cdot (4,3,2,1)$ (preserving the base octagon $\mathcal{C}$); and

$\sigma_2 = \rho_3\rho_2\rho_1\rho_3 = (1,1,1,1) \cdot (1,2,4)$ (preserving the base vertex $v$ on $\mathcal{C}$).  

The order of $\Gamma(\M)^+$ \textit{must} then be twice the number of edges in $\M$, namely $48$.
Let us assemble these and further observations in 

\begin{proposition}\label{mapstuff}
\emph{(a)}  The $3$-polytope  $\M$ is abstractly regular of type $\{8,3\}$, here realized in $\mathbb{E}^4$
 in a geometrically chiral way.
 
\emph{(b)}  The  rotation subgroup $\Gamma(\M)^+ = \langle \sigma_1, \sigma_2\rangle$ has order $48$ and presentation

\begin{equation}\label{presmpl}
\langle \sigma_1, \sigma_2\,|\,  \sigma_1^8 = \sigma_2^3 = (\sigma_1\sigma_2)^2 =
(\sigma_1^{-3}\sigma_2)^2 = 1 \rangle
\end{equation}

\emph{(c)} The full automorphism group $\Gamma(\M)$ has order $96$ and presentation

\begin{equation}\label{presm}
\langle \tau_0, \tau_1, \tau_2\,|\,  \tau_0^2 =\tau_1^2 =\tau_2^2 
= (\tau_0\tau_1)^{8} =  (\tau_1\tau_2)^{3} = (\tau_0\tau_2)^{2} = ((\tau_1\tau_0)^3 \tau_1\tau_2)^{2} =1 \rangle
\end{equation}

\end{proposition}
\noindent\textbf{Proof}. We begin with (b), where 
 it is easy to check that the  relations in (\ref{presmpl}) do hold for the matrix group 
$\langle \sigma_1, \sigma_2\rangle$. By a straightforward coset enumeration
\cite[Chapter 2]{coxmos}, we conclude from the presentation in (\ref{presmpl}) that the subgroup $\langle \sigma_1\rangle$
has the $6$ coset representatives
$$ 1, \sigma_1, \sigma_1^2,  \sigma_2\sigma_1^{-1},  \sigma_2^2\sigma_1, \sigma_2\sigma_1^{-1}\sigma_2.$$
(We abuse notation by  passing freely between the matrix group and  abstract group.) This finishes (b).

We next  note that  $\langle \sigma_1\rangle \cap \langle \sigma_2\rangle = \{1\}$, 
since $\sigma_1^j$ fixes $v$ only for $j \equiv 0 \pmod{8}$. Now we are justified in invoking
\cite[Theorem 1(c)]{schulte1}, whereby the $3$-polytope $\M$ is regular (rather than just chiral)  
if and only if  the mapping $\sigma_1 \mapsto \sigma_1^{-1}, \; \sigma_2 \mapsto \sigma_1^{2}\sigma_2$
induces  an involutory automorphism $\tau$ of $\Gamma(\M)^+$.
But the new  relations induced by applying the mapping to (\ref{presmpl}) are easily verified
formally, or even by matrices. For instance,  since $\sigma_1\sigma_2 = \sigma_2^{-1}\sigma_1^{-1}$,
we have 
$$(\sigma_1^2\sigma_2)^3 = (\sigma_1\sigma_2^{-1}\sigma_1^{-1})^3 = \sigma_1\sigma_2^{-3}\sigma_1^{-1}  = 1.$$
Thus $\M$ is abstractly regular and $\Gamma(\M)$ has order $96$. The presentation in 
(\ref{presm}) follows at once by extending $\Gamma(\M)^+$ by $\langle \tau\rangle$, then letting
$\tau_0 := \tau, \tau_1 := \tau \sigma_1, \tau_0 := \tau \sigma_1\sigma_2$.

It remains to check that our realization  is geometrically chiral.  This means that 
$\tau$ is not represented by a symmetry of $\M$ as realized in  $\mathbb{E}^4$. From the combinatorial structure, $\tau$ would have to swap vertices $1$ and $v$ while preserving the two Petrie polygons on that edge. 
This means that $\tau$ would have to act just like  $\mu_0$, that is, just like   reflection
in  a vertical line in Figure~\ref{proj1}. But $\mu_0$ does not preserve the set of $8$ edges deleted to give $\Lm$
in Figure~\ref{levi}.
\hfill$\square$

\begin{remark}
It is helpful to note that the centre of 
$\Gamma(\M)^+$ is generated by $\sigma_1^4$. 
Referring to \cite[Section 6.6]{coxmos}, we find that $\Gamma(\M)^+$ is isomorphic to the group 
$\langle -3,4 | 2\rangle$,
which in turn is an extension by $C_2$ of the binary tetrahedral group $\langle 3,3,2\rangle$.
Indeed, $a = \sigma_1^{-1} \sigma_2 \sigma_1^{-1}, b = \sigma_2 \sigma_1^4$  
satisfy $ a^3 = b^3 = (ab)^2\;  (= \zeta)$. Thus, $\langle 3,3,2\rangle \triangleleft \Gamma(\M)^+$.
\end{remark}

\section{Roli's cube -- a  chiral polytope $\R$ of type $\{8,3,3\}$}\label{rolq}
 Under the action of $G^+$ we expect to find 
$4 = 192/48$  copies of $\M$.  To understand this better, recall that there are
$12$ Petrie polygons in one chiral class, say $R$.  As with $\mathcal{C}$ and $\mathcal{C}^*$,
each polygon $\mathcal{D}$ is paired with a unique polygon $\mathcal{D}^*$ (with the disjoint set of $8$ vertices).
For each $\mathcal{D}$ there are then \textit{two} ways to remove $8$ edges so as to get a copy
of $\Lm$ and hence a copy of $\M$. Since $\M$ has six  $2$-faces like $\mathcal{C}$, we  once more find
$12\cdot2/6 = 4$ copies of $\M$.

Each Petrie polygon lies on $2$ copies of $\M$, again from the two
ways to remove $8$ edges. For example, $\mathcal{C}$ lies on both $\M$ and $(\M)\mu_0$.
(The same is true for $\mathcal{C}^*$.)

The pointwise stabilizer in $G^+$ of the base edge joining
$v = (1,1,1,1)$ and $(v)\mu_0 = (-1,1,1,1)$ must consist of
pure, unsigned even permutations of $\{2,3,4\}$. Therefore it is generated by
$$\sigma_3 := \rho_2\rho_3 = (1,1,1,1)\cdot(2,4,3).$$
It is easy to check  that $G^+ = \langle \sigma_1, \sigma_2, \sigma_3 \rangle$. 

Since three consecutive edges of a Petrie polygon lie on 
two adjacent square faces in a cubical facet of of $\Po$, it must be that every vertex
of $\R$ has the same vertex-figure as $\Po$, thus of tetrahedral type $\{3,3\}$.

We have enumerated and (implicitly) assembled the faces of a $4$-polytope $\R$,
faithfully realized in $\mathbb{E}^4$ and symmetric under the action of $G^+$. 
Let's take stock of its proper faces:

\medskip
\begin{center}
\begin{tabular}{c|c|c|c|c}
rank & stabilizer in $G^+$ &order& number of faces& type \\ \hline
$0$& $\langle  \sigma_2 , \sigma_3 \rangle$ & 12& $16$ & vertex of cube $\Po$\\ \hline
$1$& $\langle \sigma_1 \sigma_2, \sigma_3 \rangle$ & 6&$32$ & edge of $\Po$\\ \hline
$2$& $\langle \sigma_1,  \sigma_2 \sigma_3 \rangle$ &16 &$12$ & Petrie polygons of $\Po$ in one class $R$ \\ \hline
$3$& $\langle \sigma_1,  \sigma_2 \rangle$ & 48& $4$ & copy of $\M$\end{tabular}
\end{center}
It is not hard to see that our $4$-polytope $\R$ is isomorphic to Roli's cube, as  constructed in
\cite{bracho:2014aa} and as described in Section~\ref{copo}. 

\begin{theorem}\label{rconst}
\emph{(a)} The $4$-polytope $\R$ is abstractly chiral of type $\{8,3,3\}$. Its symmetry group
$\Gamma(\R) \simeq G^+ $ has order $192$ and the presentation

\begin{eqnarray}
\langle \sigma_1, \sigma_2, \sigma_3   & |& \sigma_1^8 =  \sigma_2^3 = \sigma_3^3 =  (\sigma_1 \sigma_2)^2 =  (\sigma_2 \sigma_3)^2 = 
(\sigma_1 \sigma_2 \sigma_3)^2  = 1\label{rolipresA}\\
& &\hspace*{20mm} (\sigma_1^{-3}\sigma_2)^2 = 1 \label{rolipresB}\\
& & \hspace*{20mm} (\sigma_1^{-1}\sigma_3)^4  = 1 \;\;\rangle \label{rolipresC}
\end{eqnarray}

\emph{(b)} $\R$ is faithfully realized as a geometrically chiral polytope in $\mathbb{E}^4$.

\end{theorem}
\noindent\textbf{Proof}.  The relations in (\ref{rolipresA}) are standard for chiral $4$-polytopes \cite[Theorem 1]{schulte1}; and we have seen that the relation in (\ref{rolipresB})  is a special feature of the facet $\M$.
Enumerating cosets of the subgroup $\langle \sigma_1, \sigma_2\rangle$, which still 
has order $48$, we find at most the $8$ cosets represented by
$$1, \sigma_3, \sigma_3^2,  \sigma_3^2  \sigma_1,   
 \sigma_3^2 \sigma_1^2,  \sigma_3^2 \sigma_1^2 \sigma_2,  \sigma_3^2\sigma_1^2 \sigma_2^2,  \sigma_3^2 \sigma_1^2\sigma_3 .$$
Thus  the group defined by 
(\ref{rolipresA}) and (\ref{rolipresB}) has  order at most $384$. But $G^+$, where these relations do hold, has
order $192$. We require an independent relation. In Section~\ref{mincov}, we will see why  (\ref{rolipresC}) 
is just what we need.

To show that $\R$ is abstractly chiral we must demonstrate that the mapping
$ \sigma_1 \mapsto \sigma_1^{-1},  \sigma_2 \mapsto \sigma_1^2\sigma_2, \sigma_3\mapsto\sigma_3$
 does not extend to an automorphism of $G^+$. This is easy, since
 %
 %
 \begin{equation}\label{chirpr}
 (\sigma_1\sigma_3)^4 = \zeta \;\mathrm{whereas}\;(\sigma_1^{-1} \sigma_3)^4 = 1.
 \end{equation}

Clearly, $\R$ is realized in a geometrically chiral way in
$\mathbb{E}^4$; we have already seen this for its facet $\M$.

Our concrete geometrical arguments should suffice to convince the reader that 
we really have described here a chiral $4$-polytope identical to the original
Roli's cube. A skeptic  can nail home the proof by applying  \cite[Theorem 1]{schulte1} to the group $G^+$, as generated above. 
\hfill$\square$

\bigskip
\section{Realizing the Minimal Regular Cover of $\R$}\label{mincov}

The rotation group $G^+$ for the cube has order $192$ and `standard' generators
 $\rho_0\rho_1, \rho_1\rho_2, \rho_2\rho_3$.
But for our purposes we use either of two alternate sets of generators. We already have
\begin{equation}\label{altgen1}
\sigma_1 = (\rho_0\rho_1)(\rho_2\rho_3),\; \sigma_2 = (\rho_3\rho_2)(\rho_1\rho_2)(\rho_2\rho_3),\; \sigma_3 =\rho_2\rho_3.
\end{equation}
Now we also want 
\begin{equation}\label{altgen2}
\bar{\sigma}_1 = \sigma_1^{-1},\; \bar{\sigma}_2 = \sigma_1^{2}\sigma_2,\;\bar{\sigma}_3 = \sigma_3.
\end{equation}
Recalling our the shorthand for such matrices, we have
$$ \sigma_1 =  (-1,1,1,1)\cdot(4,3,2,1);\; \sigma_2 = (1,1,1,1)\cdot(1,2,4)  ;\; \sigma_3 = (1,1,1,1)\cdot(2,4,3),
$$
and
$$\bar{\sigma}_1 = (1,1,1,-1)\cdot(1,2,3,4);\; \bar{\sigma}_2 = (-1,-1,1,1)(1,3,2); \;\bar{\sigma}_3 = (1,1,1,1)\cdot(2,4,3).
$$

We have seen that  
the group $G^+ = \langle \sigma_1, \sigma_2, \sigma_3\rangle$ (with these specified generators)
is the rotation (and full automorphism) group 
of the  chiral polytope $\R$ of type $\{8,3,3\}$. From \cite[Section 3]{schulte2} we have that the (differently generated) group 
$\bar{G^+} = \langle \bar{\sigma}_1, \bar{\sigma}_2, \bar{\sigma}_3\rangle$ is the automorphism group
for the  \textit{enantiomorphic} chiral polytope  $\bar{\R}$ .
By generating  the common group in these two ways we effectively exhibit 
right- and left-handed versions of the same polytope. 

Our geometrical realization of $\R$ began with the  base vertex $v = (1,1,1,1)$  (which
also served as base vertex for the $4$-cube $\Po$). It is crucial here that $v$ does span  the 
subspace fixed by $\sigma_1$ and $\sigma_2$. By instead taking 
$\bar{G^+}$ with base vertex $\bar{v} = (-1,1,1,1)$ fixed by  $ \bar{\sigma}_2,$ and $\bar{\sigma}_3\rangle$,
we have a faithful geometric realization of $\bar{R}$, still in $\mathbb{E}^4$, of course.

We will soon have good reason to 
 mix $G^+$ and $\bar{G^+}$ in a geometric way. Each group acts irreducibly on $\mathbb{E}^4$.
Construct the block matrices, $\kappa_j = (\sigma_j, \bar{\sigma}_j)$, $j = 1,2,3$, now acting on
 $\mathbb{E}^8$ and preserving two orthogonal subspaces of dimension $4$. 
Obviously  we may extend our notation for signed permutation matrices to the cubical group $B_8$ acting
on $\mathbb{E}^8$. Thus, taking the second copy of $\mathbb{E}^4$ to have basis $b_5, b_6, b_7, b_8$, we 
may combine our descriptions of $\sigma_j, \bar{\sigma}_j$  to get
\begin{eqnarray*}
\kappa_1 &= &(-1,1,1,1,1,1,1,-1) \cdot(4,3,2,1)(5,6,7,8),\\
\kappa_2 &= &(1,1,1,1,-1,-1,1,1) \cdot(1,2,4)(5,7,6),\\
\kappa_3 &= &\;\;\;\;\;(1,1,1,1,1,1,1,1) \cdot (2,4,3)(6,8,7).
\end{eqnarray*}

Now let  $T^+ = \langle \kappa_1, \kappa_2, \kappa_3\rangle$.
In slot-wise fashion, $\kappa_1,\kappa_2, \kappa_3$ satisfy relations like those
in (\ref{rolipresA}) and (\ref{rolipresB}).  From the proof of Theorem~\ref{rconst}, we conclude that 
$T^+$ has order $384$.  We even get a presentation for it. 

Recall that the centre of $G^+$ is generated by $\zeta = \sigma_1^4 = \bar{\sigma}_1^4$.
Thus the centre of $T^+$ has order $4$, with non-trivial elements
\begin{equation}\label{Tcen}
(\zeta,1) = (\kappa_1\kappa_3)^4,\;  (1, \zeta) = (\kappa_1^{-1}\kappa_3)^4,\mathrm{and}\;
(\zeta,\zeta) = \kappa_1^4.
\end{equation}
(This is at the heart of the proof that $\R$
is abstractly chiral.)  Looking at (\ref{chirpr}), we see that
$$ T^+/\langle (1,\zeta)\rangle \simeq \Gamma(\R),$$
and thus see the reason for the special relation in (\ref{rolipresC}).
Similarly, $ T^+/\langle (\zeta,1)\rangle \simeq \Gamma(\bar{\R})$.
Finally, we have
\begin{equation}\label{rotqu}
T^+/\langle (\zeta,\zeta)\rangle \simeq \Gamma(\Po)^+,
\end{equation}
the rotation group of the $4$-cube (isomorphic to $G^+$ generated in the customary way).

Now $T^+$ is clearly isomorphic to
the \textit{mix} $G^+ \diamondsuit\, \bar{G^+}$ 
described in \cite[Theorem 7.2]{mixa}. Guided by that  result, we 
seek an isometry $\tau_0$ of $\mathbb{E}^8$ which swaps the two orthogonal subspaces, while conjugating
each $\sigma_j$ to $\bar{\sigma}_j$. It is easy to check that
$$ \tau_0 = (1,1,1,1,1,1,1,1)\cdot (1,5)(2,6)(3,7)(4,8) $$
does the job. 
We find that  $T^+$ is the rotation subgroup of  a string C-group
 $ T = \langle \tau_0, \tau_1, \tau_2 , \tau_3\rangle$, where
$\tau_1=\tau_0 \kappa_1, \tau_2=\tau_0 \kappa_1\kappa_2, \tau_3=\tau_0 \kappa_1\kappa_2\kappa_3$. 
The corresponding directly regular
 $4$-polytope  has type $\{8,3,3\}$ and
must be the minimal regular cover of each of the chiral polytopes $\R$ and $\bar{\R}$. We consolidate 
all this in 

\bigskip
\begin{theorem}\label{mincovth} 
\emph{(a)} The group $ T = \langle \tau_0, \tau_1, \tau_2 , \tau_3\rangle$ is a string C-group
of order $768$ and with the presentation

\begin{eqnarray*} 
\langle \tau_0,\tau_1,\tau_2,\tau_3\ & |\ &  \tau_j^2= (\tau_0\tau_1)^8 = (\tau_1\tau_2)^3 =(\tau_3\tau_3)^3 = 1,\;0\leq j \leq 3,\\
& &  (\tau_0\tau_2)^2 = (\tau_0\tau_3)^2 =(\tau_1\tau_3)^2 = ((\tau_1\tau_0)^3\tau_1\tau_2)^2=1 \rangle
\end{eqnarray*}

\emph{(b)} The corresponding regular $4$-polytope $\mathcal{T}$ has type $\{8,3,3\}$ and is 
faithfully realized in $\mathbb{E}^8$,
with base vertex $(v,\bar{v}) = (1,1,1,1,-1,1,1,1,1)$.
 The polytope 
$\mathcal{T}$  is the minimal regular cover
for Roli's cube $\R$ and its enantiomorph $\bar{R}$. It  is also a  double cover of the $4$-cube $\Po$.

\emph{(c)} $\mathcal{T}  \simeq \{ \mathcal{M}\, , \, \{3,3\} \}$ is the universal regular polytope with facets $\M$ and tetrahedal vertex-figures.
\end{theorem}
\noindent\textbf{Proof}. The centre of $T$ is generated by $(\zeta,\zeta)$. It is easy to check that 
$T/\langle(\zeta,\zeta)\rangle \simeq G$, the full symmetry group of the cube; compare  (\ref{rotqu}).
In other words, the mapping $\tau_j\mapsto\rho_j,\;0\leq j \leq 3$, induces an epimorphism
$\varphi: T \rightarrow G$. Since $\tau_1,\tau_2,\tau_3$ and $\rho_1,\rho_2,\rho_3$ both satisfy the defining relations  for $\Gamma(\{3,3\}) \simeq S_4$, $\varphi$ is one-to-one on 
$\langle \tau_1,\tau_2,\tau_3 \rangle$. By the quotient criterion in \cite[2E17]{arp},
$T$ really is a string C-group. The remaining details  are routine. For background on
(c) we refer to \cite[4A]{arp}.
\hfill$\square$

\medskip

Much as in the proof, the assignment $\kappa_j\mapsto\sigma_j,\,(j=1,2,3)$,
induces an epimorphism $\varphi_R: T^+ \rightarrow G^+$. On the abstract  level, this in 
turn induces a \textit{covering} $\widetilde{\varphi}_R: \mathcal{T}^+ \rightarrow \R$, in other words,
 a rank- and adjacency-preserving
surjection of polytopes as  partially ordered sets.  The  corresponding covering of geometric polytopes
is induced by the projection 

\begin{eqnarray*}
\mathbb{E}^8 & \rightarrow & \mathbb{E}^4\\
(x,y) & \mapsto & x
\end{eqnarray*}
The projection $(x,y)\mapsto y$ likewise induces the geometrical 
covering $\widetilde{\varphi}_L: \mathcal{T}^+ \rightarrow \bar{\R}$.
Both $\widetilde{\varphi}_R$ and  $\widetilde{\varphi}_L$  are $3$-coverings, meaning here  that each 
acts isomorphically  on facets $\M$ and vertex-figures $\{3,3\}$ \cite[page 43]{arp}. Notice that  each face of
$\R$ and $\bar{\R}$ has two preimages in $\mathcal{T}$.

The polytope  $\mathcal{T}$ is also a double cover of the $4$-cube $\Po$. But there is no
natural way to embed $\Po$ in $\mathbb{E}^8$ to illustrate the geometric covering,
since $\kappa_1^4 = -1$ on any subspace of $\mathbb{E}^8$, whereas $(\rho_0\rho_1)^4 = 1$
for $\Po$.


\medskip

\section{Conclusion -  the M\"{o}bius-Kantor configuration again}

We noted earlier that $8_3$ can be `realized' as a point-line configuration in 
$\mathbb{C}^2$.  We will show this  here by first endowing $\mathbb{E}^4$ with a 
\textit{complex structure}. Thus, we want a suitable  orthogonal transformation $J$
on  $\mathbb{E}^4$ such that $J^2 = \zeta$. 
Keeping the addition, we  then  define
$$ (a+\imath b) u = au + b(uJ),\;\mathrm{for}\;  a,b\in \mathbb{R},\; u \in \mathbb{E}^4.$$
Thus $\imath u = uJ$. Over $\mathbb{C}$, $\mathbb{E}^4$ has dimension $2$.
Our choice for the matrix  $J$ is motivated by an orthogonal projection different from that 
in Figures~\ref{proj1} and \ref{levi}.

 The vectors representing the vertices labelled $0,\ldots,7$ in Figure~\ref{levi} are either  opposite 
 or  perpendicular.  Thus,   these eight points are the vertices
of a cross-polytope $\Oh = \{3,3,4\}$, one of two inscribed in $\Po$.
In \cite[Figure 4.2A]{rcp}, Coxeter gives a projection of  $\Oh$ which nicely displays
certain $2$-faces of $\Oh$.

\begin{figure}[ht]\begin{center}
\includegraphics[width=60mm]{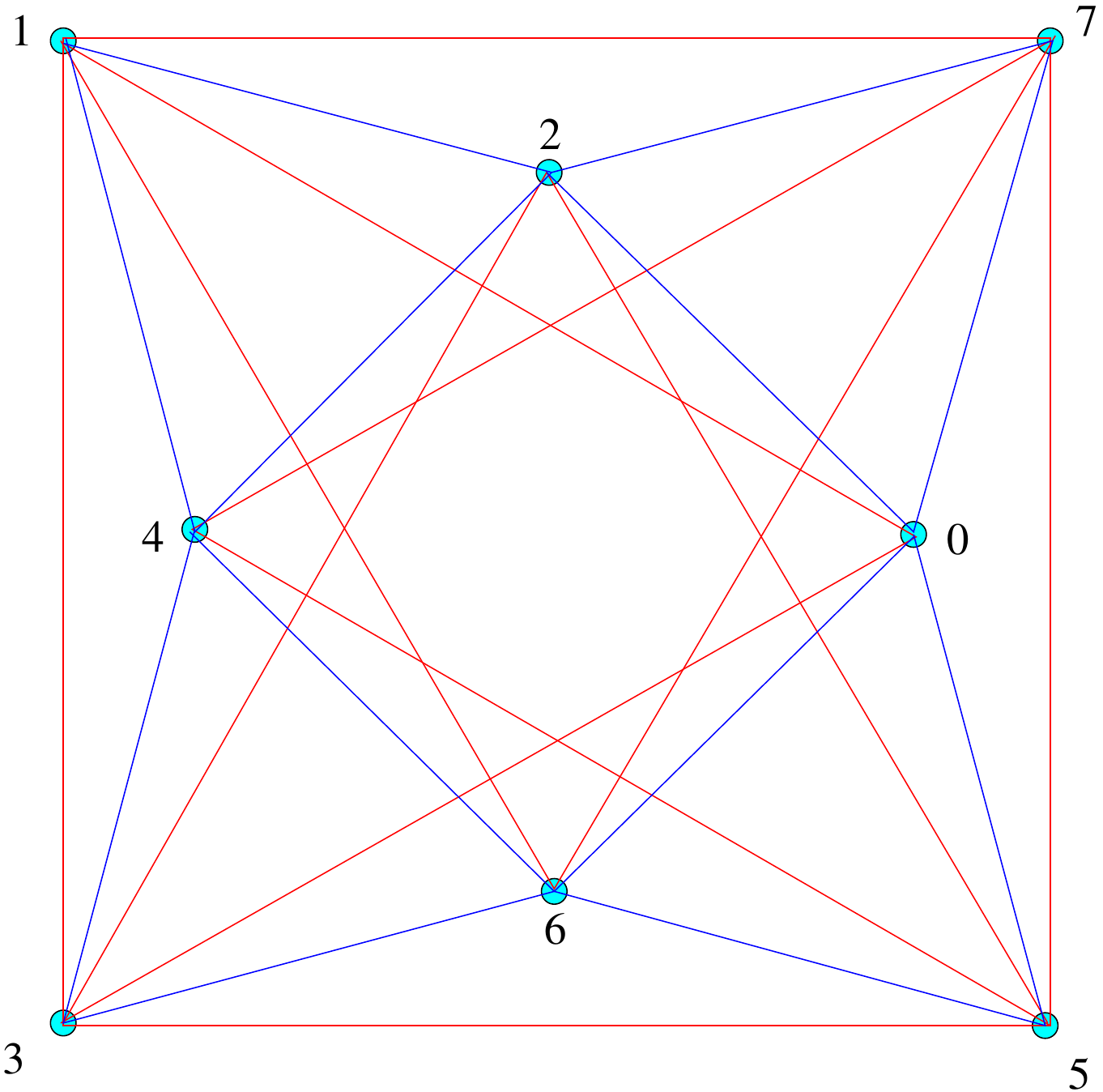}
\end{center}
\caption{Another projection of the cross-polytope $\Oh$.}
\label{cross}
\end{figure}
In Figure~\ref{cross},  each vertex of either of the two concentric squares
forms an equilateral triangle with one edge of the other square. These $8$ triangles
correspond to the unlabelled nodes in Figure~\ref{levi}, and also to the lines of 
the configuration $8_3$. (Any real triangle lies on a unique  complex line in $\mathbb{C}^2$.)
We may take the vertices in Figure~\ref{cross} to be 
$(\pm1,\pm1)$ and  $(\pm r,0), (0,\pm r)$, where $r = \sqrt{3}-1$. 

But what plane $\Lambda$ in $\mathbb{E}^4$ actually gives such a  projection?
Starting with an unknown basis $a_1, b_1$ for $\Lambda$, we  can force a lot. For example, 
edge $[2,0]$ is the projection of $(0,2,-2,0)$ and is obtained from $[7,0]$, 
the projection of $(2,2,0,0)$,   by a rotation through $60^{\circ}$. From such  details in  the geometry, 
we soon find that $\Lambda$ is uniquely determined and  get  a basis
satisfying $a_1\cdot a_1 = b_1\cdot b_1$ and $a_1\cdot b_1 = 0$. But any such  basis can still be rescaled or rotated  within $\Lambda$. Tweaking these finer details, we find it convenient  to take $a_1, b_1$ to be the first two 
rows of the matrix 

\begin{equation}\label{Lbasis}
L = \frac{1}{2\sqrt{3}}\left[\begin{array}{cccc}
        \sqrt{3}& -1 &1 & -(2+\sqrt{3})\\
          -1 &  2+\sqrt{3} & \sqrt{3}&  -1 \\
        -1 & -\sqrt{3} & 2+\sqrt{3} & 1\\
        2+\sqrt{3} & 1 & 1  & \sqrt{3}  
       \end{array}
\right]\;.
\end{equation}
The last two rows $a_2,b_2$ of $L$ give a basis for the orthogonal complement
$\Lambda^{\perp}$. 

Since we want $J$  to induce $90^{\circ}$ rotations in both $\Lambda$ and $\Lambda^{\perp}$, we
have 
\begin{equation}\label{Jmat}
J = \frac{1}{\sqrt{3}}\left[\begin{array}{rrrr}
        0&1&-1&-1\\-1&0&-1&1\\1&1&0&1\\1&-1&-1&0
       \end{array}
\right]\;.
\end{equation}
Notice that $a_1 J = b_1$ and $a_2 J = b_2$, so ,  $\{a_1, a_2\}$ is a $\mathbb{C}$-basis 
for  $\mathbb{E}^4$;   and the plane $\Lambda$ in  Figure~\ref{cross} is just $z_2=0$ in the resulting complex
coordinates. The points in the configuration $8_3$ now have these complex coordinates:

\medskip
\begin{center}
\begin{tabular}{c|c}
Label & $(z_1,z_2)$ \\ \hline
$0$ & $(  r ,  1-\imath )$  \\ \hline
$1$ & $( -1+\imath  , r  )$ \\ \hline
$2$ & $(  r\imath , -1-\imath  )$ \\ \hline
$3$ & $(  -1-\imath ,  -r\imath )$ \\ \hline
$4$ & $(  -r ,  -1+\imath )$ \\ \hline
$5$ & $(  1-\imath ,  -r )$ \\ \hline
$6$ & $(  -r\imath , 1+\imath  )$ \\ \hline
$0$ & $(  1+\imath ,  r\imath )$ 
\end{tabular}
\end{center}
(Recall that $r = \sqrt{3} - 1$.) The first coordinates do give the points displayed in Figure~\ref{cross}.
The second coordinates describe the projection onto $\Lambda^{\perp}$ ( $z_1=0$). There
  labels on the inner and outer squares  are suitably swapped.

%

A typical  line in the configuration $8_3$, like that containing points $1,6,7$, has equation
$$ r(1-\imath) z_1 + 2 z_2  = 2 r (1+\imath)\;.$$

After consulting \cite[Sections 10.6 and 11.2]{rcp}, we observe that the eight points are also the vertices of the regular complex polygon $3\{3\}3$. Its symmetry group (of unitary transformations on $\mathbb{C}^2$)
is the group $3[3]3$ with the presentation
\begin{equation}\label{333pres}
\langle \gamma_1, \gamma_2\,|\, \gamma_1^3 = 1,\; \gamma_1 \gamma_2 \gamma_1 = \gamma_2 \gamma_1 \gamma_2\rangle\;. 
\end{equation}
In fact, this group of order $24$ is isomorphic to the \textit{binary tetrahedral} group $\langle 3,3,2\rangle$.
 But in our context, we may identify it with the centralizer in $G$ of the structure matrix $J$.  A bit of computation shows that this subgroup
of $G$ is generated by 
$$ \gamma_1 = \rho_1\rho_2\rho_3\rho_2 = (1,4,2)\;\mathrm{and}\; \gamma_2 = \rho_2\rho_0\rho_1\rho_0 =
(-1,1,-1,1)\cdot (1,2,3), $$
which do satisfy the relations in (\ref{333pres}). 

 \bigskip
 \noindent\textbf{Acknowledgements}. 
I want to thank   Daniel Pellicer, both for his many geometrical ideas and
also  for generously welcoming me to the 
Centro de Ciencias Matem\'{a}ticas  at UNAM (Morelia).

\bigskip
 \bibliographystyle{siam}

\end{document}